\documentstyle{article}
\input epsf.tex

\begin{document}

\title{Minkowski-type and Alexandrov-type theorems for polyhedral herissons
\thanks{This work is partially supported by RFBR, grant 01-01-00012.}}
\author{Victor Alexandrov
\thanks{Sobolev Institute of Mathematics, Novosibirsk, Russia.
E-mail: alex@math.nsc.ru}}
\date{November 18, 2002}

\maketitle
\begin{abstract}
Classical H.~Minkowski theorems on existence and uniqueness of convex
polyhedra with prescribed directions and areas of faces as well as the well-known
generalization of H.~Minkowski uniqueness theorem due to A.~D.~Alexandrov
are extended to a class of nonconvex polyhedra which are called polyhedral
herissons and may be described as polyhedra with injective spherical image.

Key words: convex polyhedron, polyhedral surface, polyhedral hedgehog,
equipment, virtual polytope, polygon, Cauchy lemma, open mapping.

2000 Mathematics Subject Classification: 52B10, 52C25, 52B70, 52A38, 52A15.
\end{abstract}

\section{Introduction}
H.~Minkowski proved that a convex polyhedron is uniquely (up to translation)
determined by the areas of its faces and the unit outward normal vectors to the faces
(see, for example, \cite{Al} or \cite{Mi}).
Moreover, H.~Minkowski found natural and easily verified conditions on
a finite set of unit vectors and positive numbers implying that there exists
a convex polytope $P$ such that these vectors are the outward normal vectors to
the faces of $P$ and these numbers equal the areas of the faces of $P$.

We extend these results (as well as a generalization of the Minkowski uniqueness
theorem which is due to A.~D.~Alexandrov) to a class of nonconvex polyhedral
surfaces.
These surfaces may even admit self-intersections.
The crucial point is that their spherical mapping is injective.
They are called polyhedral herissons (``herisson'' is a French word for ``hedgehog'')
and were previously studied, for example, in papers by
J.~J.~Stoker \cite{St},
L.~Rodrigues and H.~Rosenberg \cite{RR},
P.~Roitman \cite{Ro}
and G.~Yu.~Panina \cite{Pa}.

\section{H.~Minkowski and A.~D.~Alexandrov theorems for convex polyhedra}

In this section we make some notions more precise and revise the classical
H.~Minkowski and A.~D.~Alexandrov theorems on existence and uniqueness of convex
polyhedra with prescribed directions and areas of faces.

{\bf 2.1. Definitions.}  Let $L$ be an affine subspace of ${\bf R}^d$, $d\geq 2$
(the case $L={\bf R}^d$ is not excluded).
The convex hull of a finite set of points in $L$ is said to be a {\it convex polytope}
provided that it has interior points in $L$.
A convex polytope in a two-dimensional space $L$ is called a {\it convex polygon}.

Given a convex polytope $P$ in ${\bf R}^d$  and a unit vector
$n\in{\bf R}^d$, find the least real number $h$ such that, for each
$h'>h$, the hyperplane $(x,n)=h'$ does not intersect $P$.
Here and below $(x,n)$ stands for the scalar product of vectors $x$ and $n$.
The intersection of $P$ and the above-constructed hyperplane
$(x,n)=h$ is called a {\it face} of $P$ with outward normal vector $n$,
while the hyperplane $(x,n)=h$ itself is called a {\it supporting hyperplane}
to the face.
The dimension of the affine hull of a face is said to be the
{\it dimension of the face}.
Note that, under the above definition of a face, the dimension of a face
can take any integer value from 0 to $d-1$ inclusively.

Let $P_1$ and $P_2$ be convex polytopes in $R^d$, $d\geq 2$, and let $Q_1$ be a
$k$-dimensional face of $P_1$, $0\leq k\leq d-1$.
Denote by $n$ the outward unit vector to $Q_1$.
Let $Q_2$ be a face of $P_2$ with outward normal vector $n$.
Then $Q_2$ is called {\it parallel} to $Q_1$.

Let $L_j$ be an affine subspace of ${\bf R}^d$, $j=1,2$, and let $Q_j$ be
a convex polytope in $L_j$.
Say that the polytope $Q_1$ {\it can be put inside} $Q_2$ by a translation if
there exists a translation $T:{\bf R}^d \to{\bf R}^d$ such that $T(Q_1)$
is a subset of $Q_2$ but does not coincide with $Q_2$.

{\bf 2.2. Theorem} (A.~D.~Alexandrov \cite{Al}).
{\it Let $d=3$ and let $P_1, P_2$ be convex polytopes in $R^d$.
Suppose that, given any two parallel faces of
$P_1$ and $P_2$ such that at least one of the faces is of dimension $d-1$,
neither of the faces can be put inside the other face by a translation.
Then $P_1$ and $P_2$ are congruent and parallel to each other
(i.e. can be obtained from one another by a translation).}

{\bf 2.3. Remarks.}
The assertion of Theorem 2.2 is trivial for $d=2$ and, as is observed in \cite{Al},
is false for $d=4$.
As a counterexample, one can consider
a cube with edges of length 2 in ${\bf R}^4$ and a rectangular
parallelepiped with edges of lengths 1, 1, 3, and 3
parallel to the edges of the cube.
It is clear that none of the 3-dimensional faces of the cube can be put inside
the parallel face of the parallelepiped by a translation as well as none of
the 3-dimensional faces of the parallelepiped can be put inside the parallel
face of the cube by a translation.
Nevertheless it is clear that the cube and the parallelepiped are not congruent.

For $d>4$, a similar counterexample can be obtained in
${\bf R}^d={\bf R}^4\times{\bf R}^{d-4}$
by considering the Cartesian product of the above-constructed 4-dimensional
cube and parallelepiped with $d-4$ unit segments.

Note that, to apply Theorem 2.2, we should have a positive answer to the
question ``Is it possible to put the face with outward unit normal vector $n$
of one of the polytopes inside the parallel face of the other polytope?''
for only finitely many vectors $n$ (namely, only for the outward unit normal vectors
to the 2-dimensional faces).
The next theorem shows that a positive answer to the above question for all
vectors  $n\in{\bf R}^d$ implies the conclusion of Theorem 2.2 in all dimensions.
The problem remains open whether it is possible to find an analog of Theorem 2.4
where a positive answer to the above question is requested for only finitely many
vectors and how many vectors should be tested.

{\bf 2.4. Theorem.}
{\it Let $d\geq 2$ and let $P_1, P_2$ be convex polytopes in ${\bf R}^d$.
Suppose that there is no unit vector $n\in{\bf R}^d$ such that the face
of $P_1$ or $P_2$ with outward normal vector $n$ can be put
inside the parallel face of the other of the polytopes $P_1$ and $P_2$ by a
translation.
Then $P_1$ and $P_2$ are congruent and parallel to each other.}

{\bf 2.5. Proof.}
Let $v$ be a 0-dimensional face of $P_1$ and let $n$ be the unit outward normal vector
to one of the supporting hyperplanes of $P_1$ passing through $v$.
If the face $w$ of $P_2$ with outward normal vector $n$ has nonzero dimension,
we can put $v$ inside $w$ by a parallel translation because $v$ is a point.
This contradicts the conditions of Theorem 2.4.
Hence the dimension of $w$ equals zero.

Let $T:{\bf R}^d\to{\bf R}^d$ be a translation such that $v=T(w)$.
As far as $n$ is arbitrary in the above reasoning, it follows that
$P_1$ and $T(P_2)$ coincide in a neighborhood of $v$.
In particular, for each 1-dimensional face $V$ of $P_1$ that is incident to $v$,
a small neighborhood of $v$ in $V$ is contained in a 1-dimensional face $W$ of $T(P_2)$.
If the lengths of $V$ and $W$ were different,
we could put the smaller segment into the bigger one by a translation
in contradiction to the conditions of Theorem 2.4.
Hence, 1-dimensional faces $V$ and $W$ coincide and $P_1$ and $T(P_2)$
have one more common 0-dimensional face besides $v$.
Denote it by $\widetilde v$.

Examining as before all supporting planes passing through $\widetilde v$,
we conclude that $P_1$ and $T(P_2)$ coincide in a neighborhood of $\widetilde v$.
Moreover, each 1-dimensional face of $P_1$ that is incident to $\widetilde v$
coincides with a 1-dimensional face of $T(P_2)$ and vice versa.
Thus there exists one more common 0-dimensional face of $P_1$ and $T(P_2)$
besides $v$ and $\widetilde v$.

After a finite number of iterations, this travel along 0-dimensional faces of
the polytopes shows that the set of 0-dimensional faces of $P_1$ coincides with
the set of 0-dimensional faces of $T(P_2)$.
Therefore $P_1$ and $T(P_2)$ coincide and $P_1$ and $P_2$ are congruent and
parallel to each other.
$\framebox{\phantom{o}}$

{\bf 2.6. Theorem} (H.~Minkowski uniqueness theorem \cite{Al}, \cite{Mi}).
{\it Let $d\geq 2$ and let two convex polytopes in ${\bf R}^d$
be such that, for every $(d-1)$-dimensional face of each of the polytopes,
the parallel face of the other polytope has the same $(d-1)$-dimensional volume.
Then the polytopes are congruent and parallel to each other.}

{\bf 2.7. Remarks.}
The assertion of Theorem 2.6 is trivial for $d=2$ and is an immediate consequence of Theorem 2.2
for $d=3$.
Nevertheless there are proofs of Theorem 2.6 that are independent of Theorem 2.2
and valid for all $d\geq 3$.

Recall the following well-known fact:
If $d\geq 2$, $P$ is a convex polytope in ${\bf R}^d$, $n_1,\dots ,n_m$ are
the unit outward normal vectors to $(d-1)$-dimensional faces of $P$, and
$f_1,\dots ,f_m$ are the $(d-1)$-dimensional volumes of the corresponding faces,
then $\sum_{j=1}^{m} f_jn_j=0$.

{\bf 2.8. Theorem} (H.~Minkowski existence theorem \cite{Al}, \cite{Mi}).
{\it Let $d\geq 2$ and let $n_1,\dots ,n_m$ be unit vectors in ${\bf R}^d$
which do not lie in a closed half-space bounded by a hyperplane passing through the origin.
Let $f_1,\dots ,f_m$ be positive real numbers such that
$\sum_{j=1}^{m} f_jn_j=0$.
Then there exists a convex polytope $P$ in  ${\bf R}^d$ such that the vectors
$n_1,\dots ,n_m$ (and only they) are the unit outward normal vectors to
the $(d-1)$-dimensional faces of $P$ and the $(d-1)$-dimensional volumes of the corresponding
faces are equal to $f_1,\dots ,f_m$.}

In this paper we prove some analogs of Theorems 2.2, 2.6, and 2.8
for a class of nonconvex polyhedra.

\section{Polyhedral herissons}

In this section we introduce nonconvex polyhedral surfaces in ${\bf R}^3$
to which we are going to extend Theorems 2.2, 2.6, and 2.8.

{\bf 3.1. Definitions.}
The image of a 2-dimensional sphere-homeomorphic cell complex under a continuous
mapping into ${\bf R}^3$ is said to be a {\it polyhedral surface} if
the image of each cell is a convex $k$-dimensional polytope in
${\bf R}^3$, $k=0,1,2$.
In this case, the image of a $k$-dimensional cell is called a {\it face}, for $k=2$;
an {\it edge}, for $k=1$; and a {\it vertex} of the polyhedral surface, for $k=0$.

{\bf 3.2. Remark.}
A polyhedral surface in the sense of Definition 3.1 need not be the boundary of
a convex polyhedron, may have rather complicated self-intersections, but
is an orientable manifold.

{\bf 3.3. Definitions.}
A polyhedral surface $P$ in ${\bf R}^3$ is called a {\it polyhedral herisson} if
each face $Q_j$ of $P$ is equipped with a unit normal vector $n_j$ in such
a way that, if faces $Q_j$ and $Q_k$ share an edge, then $n_j + n_k \neq 0$ and
thus there is only one shortest geodesic path on the unit sphere
${\bf S}^2\subset{\bf R}^3$ joining the endpoints of $n_j$ and $n_k$.
Denote the shortest geodesic path by $l_{jk}$ and suppose in addition that,
for arbitrary indices $j$ and $k$ corresponding to adjacent faces, any two curves
$l_{jk}$ either do not intersect or have one endpoint in common.
Finally, we suppose that the union of all curves $l_{jk}$ gives rise to a partition
of the sphere into convex polygons.

The set of vectors $n_1,\dots ,n_m$ involved in the above definition is called an
{\it equipment} of the herisson.
The polyhedral herisson is denoted by $(P; n_1,\dots ,n_m)$.
The convex spherical polygon of the above-constructed partition of the sphere
${\bf S}^2$ is called the {\it spherical image} of the corresponding vertex.

{\bf 3.4. Example.}
The boundary of every convex polytope in ${\bf R}^3$ equipped with, say,
the unit outward normal vectors to the faces is a polyhedral herisson.

$$
\epsfxsize=13cm
\epsfbox{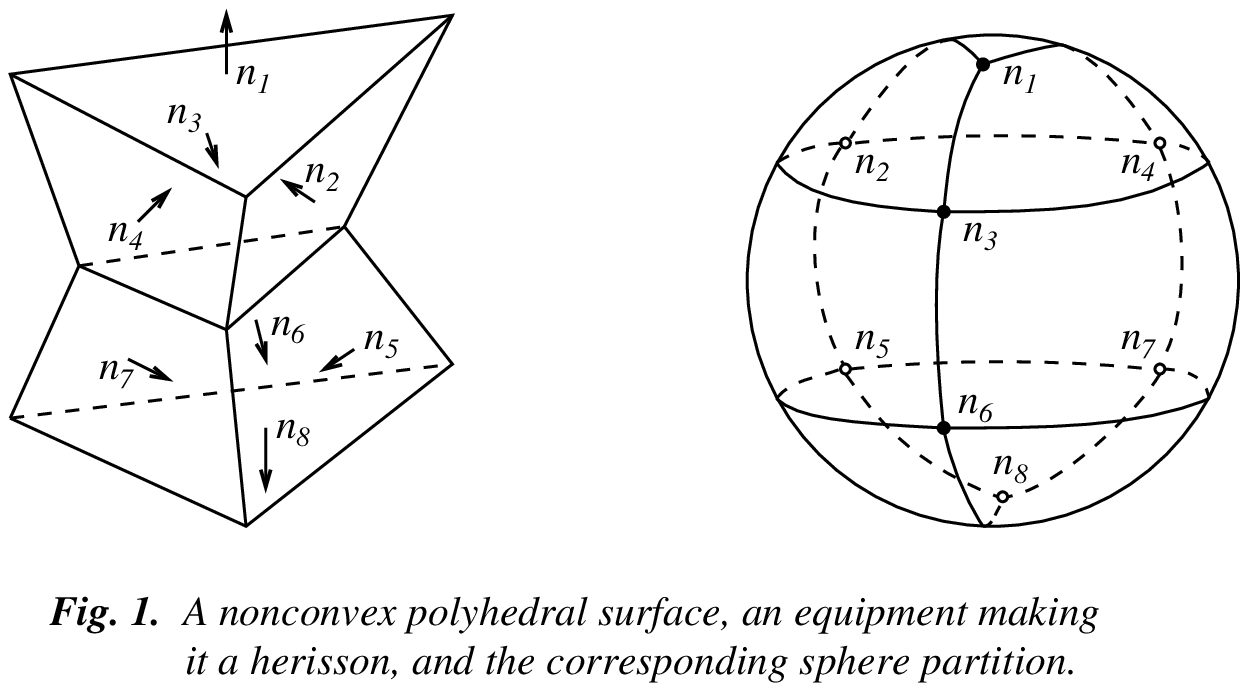}
$$

{\bf 3.5. Example.}
Consider a regular tetrahedron $\Delta$ in ${\bf R}^3$.
Cut a smaller tetrahedron from $\Delta$ by a plane $\tau$ which is parallel
to one of the faces of $\Delta$.
Consider the union of the truncated tetrahedron and its image under the
reflection in $\tau$.
Denote by $P$ the boundary of the union (Fig.~1).
Equip the faces of $P$ which are parallel to $\tau$ with the
unit outward normal vectors and equip all the rest faces of $P$ with
the unit inward normal vectors.
This equipment makes $P$ a polyhedral herisson.
The corresponding sphere partition is shown on the right-hand side of Fig.~1.

{\bf 3.6. Remarks.}
Despite
unhandiness
of the definition, polyhedral herissons can be
considered as a quite natural generalization of convex polytopes.

To clarify this point, note that every convex $C^\infty$-surface
$S$ in ${\bf R}^d$ which is homeomorphic to the sphere ${\bf S}^{d-1}$
can be considered as the envelope of its supporting hyperplanes
$(x,n)=h(n)$, $n\in{\bf S}^{d-1}$ (Fig.~2a).
Here $h$ is a $C^\infty$-function on the sphere ${\bf S}^{d-1}$.
Given a constant number $t$, one can expect that the $C^\infty$-function
${\bf S}^{d-1}\ni n\mapsto h(n)+t$ is a support function of some surface $S_t$
that is called parallel to $S$.
Both smoothness and the shape of $S_t$ strongly depend on the sign of $t$, i.e.
on  whether we build the parallel surface outside or inside of $S$ (Fig.~2b and 2c).
It may occur that the envelope of a $C^\infty$-smooth family of hyperplanes
$(x,n)=h(n)+t$, $n\in{\bf S}^{d-1}$, is neither convex nor regular hyperplane.

$$
\epsfbox{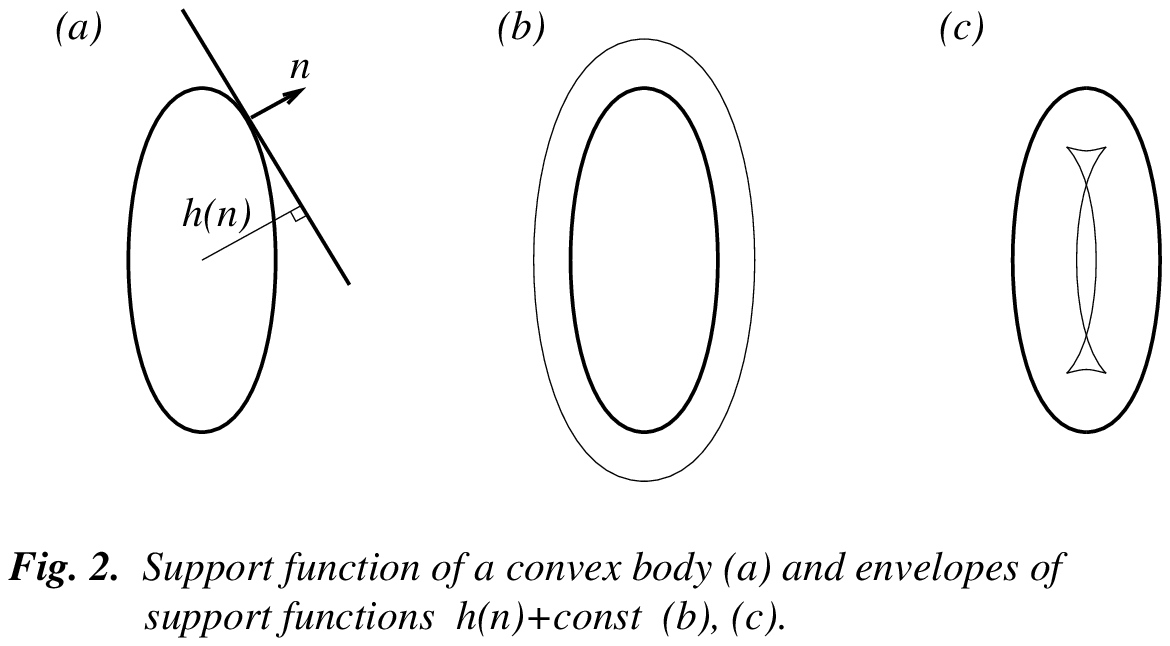}
$$

In \cite{LLR}, R.~Langevin, G.~Levitt, and H.~Rosenberg began the study of such
not necessarily convex or regular hypersurfaces which admit a
$C^\infty$-parametrization by a unit sphere.
In some sense these are hypersurfaces with one-to-one Gauss mapping
that justifies the title ``herisson'' or hedgehog: in each direction a
unique quill is sticked out.
In \cite{LLR} examples are given which show how herissons appear in the theory
of minimal surfaces, algebraic hypersurfaces, singular wave fronts, etc.
In \cite{Ma1}--\cite{Ma8} it is shown that the notion of a herisson is
useful for studying various problems of convex geometry.

Strictly speaking, the notion of a polyhedral herisson is introduced in  \cite{RR}
where several theorems are proven which are similar to the famous
Cauchy theorem on unique determination of a convex polytope by the intrinsic
metric of its boundary \cite{Al}, \cite{Ca}.
G.~Yu.~Panina established that a polyhedral herisson is a special
case of a more general and more algebraic object called a virtual polytope.
The latter can be treated as the Minkowski difference of two convex polytopes.
It appeared naturally in the theory of polytope algebras
(see \cite{Mc}, \cite{Mo}, and \cite{PKh}).
Another generalization of the notion of a polyhedral herisson (which is more
geometrical in comparison with the notion of a virtual polytope) is given
in the Ph.D. thesis of P.~Roitman \cite{Ro}.

{\bf 3.7. Definitions.}
Let $(P; n_1,\dots, n_m)$ be a polyhedral herisson in ${\bf R}^3$
and let a set of unit normal vectors $(\nu_1,\dots,\nu_m)$
to faces of the polyhedral surface $P$ be an orientation of $P$.
By definition, put $\varepsilon_j=(n_j,\nu_j)$ for each $j=1,\dots, m$.
A face $G_j$ is called {\it positive} if $\varepsilon_j=+1$ and
{\it negative} if $\varepsilon_j=-1$.
The product of $\varepsilon_j$ with a nonnegative number that is equal to the area
of the convex polygon $G_j$ is called the {\it oriented area} of the face $G_j$.

Given two polyhedral herissons in ${\bf R}^3$, their faces endowed with
a common normal vector  $n_j$ are called {\it parallel}.
The herissons are called {\it parallel and of the same orientation}
if the following conditions are fulfilled:
(i) their equipments coincide with each other;
(ii) their equipments generate the same partitions of the sphere;
(iii) every pair of their parallel faces consists either of positive or of negative faces.

We say that two polyhedral herissons are {\it congruent and parallel} to each other if
one of them (together with its equipment) is obtained from the other by a translation.

{\bf 3.8. Remarks.}
In Definition 3.7, conditions (i)--(iii) are independent.

Indeed, convex polyhedral surfaces on Fig.~3 equipped with the unit
outward normal vectors to the faces give us an example of polyhedral herissons with
a common equipment but different sphere partitions.
Hence (i) does not imply (ii).

Given a polyhedral herisson, alter the orientation of the corresponding polyhedral
surface without altering the equipment.
The initial and modified herissons possess properties (i) and (ii) but
do not possess property (iii).

$$
\epsfbox{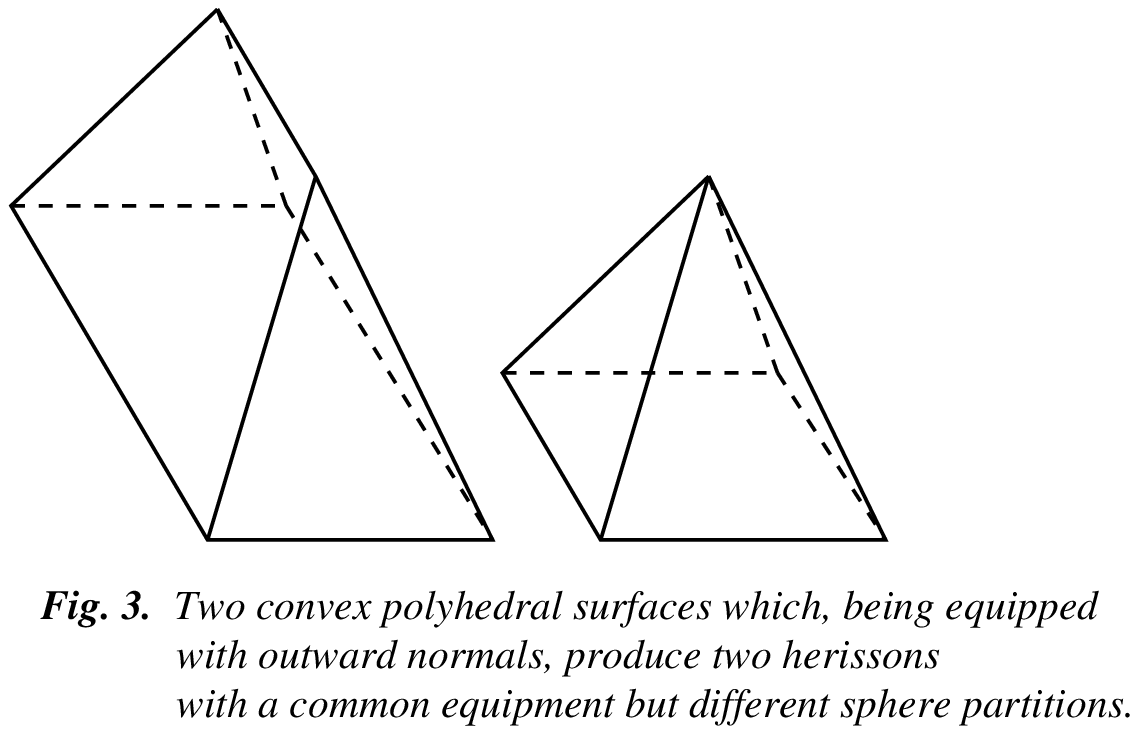}
$$

\section{Uniqueness theorems for polyhedral herissons}

{\bf 4.1. Theorem.}
{\it  Given two polyhedral herissons in ${\bf R}^3$ which are parallel and of the
same orientation, suppose  that they have no two parallel faces
such that one of the faces can be put inside the other face by a translation.
Then the herissons are congruent and parallel to each other.}

{\bf 4.2. Remarks.}
Theorem 4.1 may be considered as an analog of Theorem 2.2
though does not include it as a particular case (Theorem 2.2
works even if the dimension of one of the parallel faces is 0 or 1).
Moreover, the proof of Theorem 4.1 presented below is similar
to that of Theorem 2.2 given in  \cite{Al}.
Both proofs are based on Lemmas 4.3 and 4.4 which are given below
without proofs (the reader can find the proofs in \cite{Al} and \cite{St}).

{\bf 4.3. Lemma} (A.~Cauchy \cite[Chap.~II, \S1]{Al}).
{\it Let a sphere-homeomorphic cell complex be given such that no 2-cell
is bounded by two edges.
Let the edges of the complex be arbitrarily labeled with the numbers $+1$, $0$, or $-1$
and let the index $j$ be assigned to each vertex of the complex that is equal to the
number of changes of signs of edges from $+1$ to $-1$ and from $-1$ to $+1$
counted in making one circuit around the vertex (edges labeled with
$0$ are ignored).
The conclusion is that either all the edges of the complex are labeled with $0$ or
there is a vertex which is incident to an edge labeled with $+1$ or $-1$
and for which $j\leq 2$.}

An example is given in Fig.~4 where the edges of a cell complex isomorphic to
the boundary of a tetrahedron are labeled with the numbers $+1$, $0$, or $-1$
and the index $j$ is given for each vertex.

$$
\epsfbox{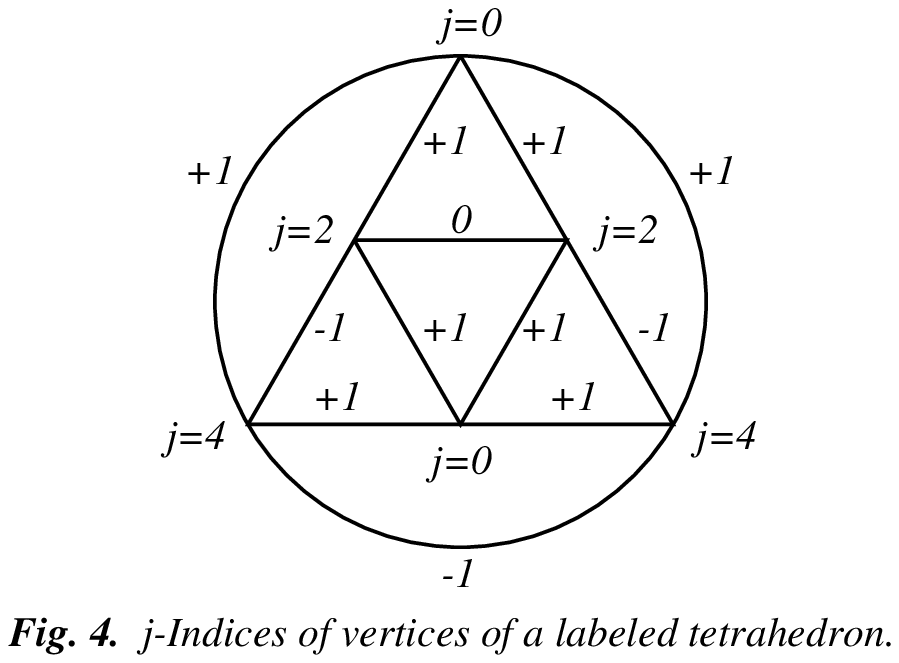}
$$

{\bf 4.4. Lemma} (A.~D.~Alexandrov \cite[Chap.~VI, \S1]{Al}).
{\it Let two convex polygons be given such that neither of them can be put
inside the other by a translation.
Let the vertices and edges of the polygons be labeled with the numbers $+1$, $0$, or $-1$
according to the following rules:

{\rm(i)}~if a vertex of one of the polygons is such that, for every outward normal vector to
the polygon at the point, the parallel face of the other polygon is a vertex,
then both vertices are labeled with $0$;

{\rm(ii)}~if a vertex of one of the polygons is such that, for some outward normal vector to
the polygon at the point, the parallel face of the other polygon is an edge,
then the vertex is labeled with $-1$ and the edge is labeled with $+1$;

{\rm(iii)}~if an edge of one of the polygons is such that the parallel face of the
other  polygon is a vertex,
then the edge is labeled with $+1$ and the vertex is labeled with $-1$;

{\rm(iv)}~finally, if an edge  of one of the polygons is such that the parallel face of the
other polygon is an edge, then the longer edge is labeled with $+1$ and the shorter
edge is labeled with $-1$; if the edges in question have equal lengths, both
of them are labeled with $0$.

Let each of the polygons be assign with the index $i$ that is equal to the number of
changes of signs of edges and vertices from $+1$ to $-1$ and from $-1$ to $+1$
counted in making one circuit around the polygon (edges and vertices labeled with
$0$ are ignored).
The conclusion is that either the polygons are congruent and parallel (and, hence,
$i=0$) or $i\geq 4$ for each of the polygons.}

An example of calculation of the index $i$ is shown in Fig.~5.

$$
\epsfbox{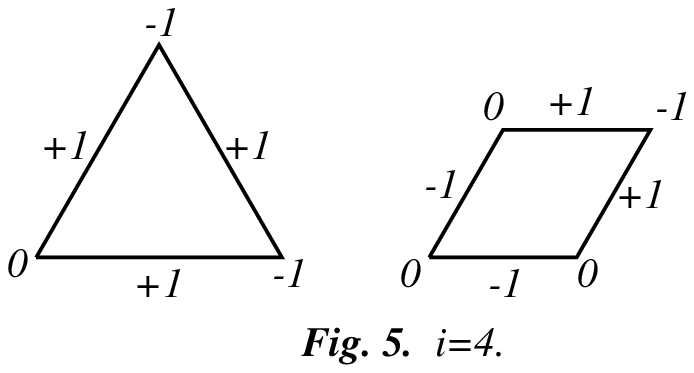}
$$

{\bf 4.5. Proof} of Theorem 4.1.
Label the edges of the two given polyhedral herissons with $+1$, $0$, or $-1$ by applying
the rules described in Lemma 4.4 to each pair of their parallel faces.
Since the herissons are parallel and of the same orientation,
it follows that parallel faces have mutually parallel edges.
Therefore each vertex is labeled with $0$.
Lemma 4.4 implies that either parallel faces are congruent and thus all
their edges are labeled with $0$ or some of their edges are labeled with
$+1$ or $-1$ and in the latter case $i\geq 4$ for these faces.

Transfer the above-constructed distribution of the numbers $+1$, $0$, and $-1$
from the edges of the polyhedral surface $Q$ of one of the herissons to the
cell complex $K$ whose image is $Q$.
Then for every 2-cell of $K$ either all 1-cells bounding the 2-cell
are labeled with $0$ or some of the 1-cells are labeled with $+1$ or $-1$
and in the latter case $i\geq 4$.

By assumption, $K$ is sphere homeomorphic and thus there exists a dual
cell complex $C$ for it.
We can treat $C$ as follows: In each face of $Q$ choose a point and
join two such points by a simple curve on $Q$ if and only if
the corresponding 2-cells share an edge on $Q$. Here we assume
that any two of these curves can intersect each other at endpoints only.
A set of 1-cells in $C$ bounds a 2-cell of $C$ if and only if
the corresponding 1-cells of $K$ are incident to a vertex and
no other 1-cell of $K$ is incident to the vertex.
All we need from the definition of the dual graph is that the edges
of $K$ and $C$ are in one-to-one correspondence.
This makes it possible to transfer the above-constructed distribution
of the numbers  $+1$, $0$, and $-1$ from the edges of $K$ to the edges of $C$.

Note that no 2-cell of $C$ is bounded by only two 1-cells of $C$.
Otherwise there would be a 0-cell $v$ of $K$ incident to only two 1-cells of $K$.
In this case there would be a vertex on $Q$ incident to only two edges and thus
to only two faces.
Denote the edges by $e_1$ and $e_2$ and the faces by $f_1$ and $f_2$.
If $e_1$ and $e_2$ do not lie on one straight line then $f_1$ and $f_2$
lie on one plane.
It follows that the corresponding vectors $n_1$ and $n_2$ of the equipment
of the herisson either coincide or are opposite.
Both possibilities contradict with Definition 3.3 of a polyhedral herisson.
If $e_1$ and $e_2$ lie on one straight line then we can modify $K$ in such a way as
to exclude $v$ and replace the two edges $e_1$ and $e_2$ by a single edge
which joins together the vertices of $e_1$ and $e_2$ that are different from $v$.
Obviously, $Q$ can be treated as a piecewise affine image of the modified complex.
On the other hand, the complex $C$ corresponding to the modified complex $K$
contains no 2-cell which is bounded by only two 1-cells.

Note that the index $j$ of a vertex of $K$ is equal to the index $i$ of the
corresponding 2-cell of $K$.
Comparing the restrictions imposed on the indices $i$ and $j$ by Lemmas 4.3 and 4.4
we conclude that both indices identically equal zero.
Moreover, we conclude that every pair of parallel faces is a pair of
congruent and parallel convex polygons.

Now prove that polyhedral herissons satisfying the
conditions of Theorem 4.1 can be obtained (together with the equipments)
from one another by a translation.
To this end superpose some two parallel faces of the herissons by a translation
in such a way that the two corresponding vectors of the equipments coincide.
Denote the coincident faces by $Q_1$.
Let $Q_2$ be a face of the first herisson which shares a common edge with $Q_1$.
The vector of the equipment attached to $Q_2$ defines a point on the sphere
${\bf S}^2$ which is joint with the vector of the equipment attached to $Q_1$.
By definition, the equipments of the herissons generate the same sphere partition.
Hence there is a face $\widetilde Q_2$ of the second herisson that is parallel
to $Q_2$ and shares an edge with $Q_1$.
Obviously, the dihedral angle between $Q_1$ and $Q_2$ is equal to the
dihedral angle between $Q_1$ and $\widetilde Q_2$.
Besides, we already know that the parallel faces $Q_2$ and $\widetilde Q_2$
can be superposed by a translation.
Thus $Q_2$ and $\widetilde Q_2$ coincide;
moreover, the attached vectors of the equipments also coincide.

If we continue this process of checking that next adjacent faces automatically
coincide, we see that, as soon as we have superposed the faces $Q_1$
of the herissons, we have superposed the herissons in total (including the
equipments).
$\framebox{\phantom{o}}$

{\bf 4.6. Theorem.}
{\it  Given two polyhedral herissons in ${\bf R}^3$ which are parallel and of the
same orientation, suppose that their parallel faces
have equal oriented areas.
Then the herissons are congruent and parallel to each other.}

{\bf 4.7. Remarks.}
Theorem 4.6 is similar to Theorem 2.6 but is valid for nonconvex polyhedral surfaces.

Theorem 4.6 is an immediate consequence of Theorem 4.1.

In the proof of Theorem 4.1 (and thus of Theorem 4.6) we do not use
that the equipment of a herisson generates a sphere partition into
just convex polygons (see Definition 3.3).

\section{Existence theorems for polyhedral herissons}

In this section we prove the main result of the paper.
Namely, we prove Theorem 5.7 that generalizes the Minkowski existence theorem
for convex polytopes with prescribed directions of faces and areas of faces
to the case of polyhedral herissons.
We begin with auxiliary results.

{\bf 5.1. Lemma.}
{\it Let $(P; n_1,\dots ,n_m)$ be a polyhedral herisson in ${\bf R}^3$,
let $(\nu_1,\dots,\nu_m)$ be an orientation of the polyhedral surface $P$,
and let $f_1,\dots,f_m$  be the oriented areas of the faces of the herisson.
Then $\sum_{j=1}^{m}f_jn_j=0$.}

{\bf 5.2.} The {\bf proof} is obtained by the direct calculation
$$
\sum_{j=1}^{m}f_jn_j=\sum_{j=1}^{m}|f_j|\varepsilon_jn_j=
\sum_{j=1}^{m}|f_j|\nu_j=0,
$$
where we use the notations of Definition 3.7 and the last equality
is written according to the well-known observation due to Minkowski
which was mentioned in Remark 2.7.
$\framebox{\phantom{o}}$

{\bf 5.3. Lemma.}
{\it Let $P_0$ be a convex polygon on the Euclidean plane such that
no two sides of $P_0$ lie on parallel straight lines.
Let $\cal S$ be a family of convex polygons $P$ satisfying the following conditions:

{\rm(a)} for every side of every $P\in\cal S$ there is a parallel side of $P_0$;

{\rm(b)} for every $P\in\cal S$ and every side of $P_0$ there is a parallel side of $P$;

{\rm(c)} the areas of all polygons $P\in\cal S$ are uniformly bounded from above.

Then the perimeters of all polygons $P\in\cal S$ are also uniformly bounded from above.}

{\bf 5.4. Remark.}
The assertion of Lemma 5.2 is not true without the assumption that $P_0$ has
no parallel sides.
To make sure of this, we can consider the square with unit-length sides
as $P_0$ and the sequence  of rectangles $P_k$ with sides parallel to the sides of
$P_0$ and side lengths $k$ and $1/k$, $k\in{\bf N}$, as $\cal S$.

{\bf 5.5.} The {\bf proof} of Lemma 5.3 is obtained by reductio ad absurdum.
Suppose the perimeters of polygons of $\cal S$ are not bounded from above.
Conditions (a) and (b) imply that every polygon $P\in\cal S$
has the same number of sides as $P_0$.
Therefore there is a sequence of polygons $P_k$ in $\cal S$ such that
the length $l_k$ of a maximal side of $P_k$ tends to infinity as $k\to\infty$.
Denote by $\alpha$ the least angle between straight lines containing two
different sides of $P_0$.
Obviously,  $\alpha >0$, since $P_0$ has no parallel sides.
For each $k\in{\bf N}$, construct an isosceles triangle $\Delta_k$ such that
the side of $P_k$ of length $l_k$ is the base of $\Delta_k$;
the angles attached to the base of $\Delta_k$ equal $\alpha$;
$\Delta_k$ and $P_k$ lie on the same side from their common side of length
$l_k$ (Fig.~6).

$$
\epsfbox{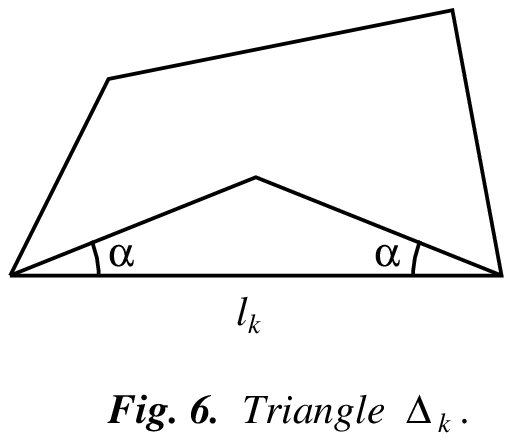}
$$

Because of the choice of $\alpha$, $\Delta_k$ is contained in $P_k$.
The latter fact leads to a contradiction.
Namely, on the one hand, the area of  $\Delta_k$ equals
$\frac{1}{4}l_k^2\sin^2\alpha$ and thus tends to infinity as $k\to\infty$.
On the other hand, the area is less than or equal to the area of $P_k$ and hence
is uniformly bounded from above.
This contradiction concludes the proof.
$\framebox{\phantom{o}}$

{\bf 5.6. Definition.}
A polyhedral herisson is said to be in a {\it general position} if
no three vectors of its equipment are coplanar.

{\bf 5.7. Theorem.}
{\it Let $(P; n_1,\dots,n_m)$ be a polyhedral herisson in ${\bf R}^3$ in a general
position.
Fix an orientation of the polyhedral surface $P$ and denote the oriented areas of
the faces of $(P; n_1,\dots,n_m)$ by $f_1,\dots,f_m$ (in particular, it is assumed that
none of the numbers $f_1,\dots,f_m$ equals zero).
Suppose $g_1,\dots,g_m$ are real numbers such that (1) $f_j\cdot g_j>0$ for every
$j=1,\dots, m$ and (2) $\sum_{j=1}^{m}g_jn_j=0$.
Then there exists a polyhedral herisson in ${\bf R}^3$ which is parallel to and
of the same orientation as $(P; n_1,\dots,n_m)$
and is such that $g_1,\dots,g_m$ are the oriented areas of the faces.}

{\bf 5.8. Remark.}
Theorem 5.7 is similar to Theorem 2.8 but, even in the 3-space,
the former does not contain the latter as a particular case.
The main difference is that in Theorem 2.8 we start from abstract sets of vectors
$n_1,\dots,n_m$ and real numbers $f_1,\dots,f_m$ which satisfy some relation.
In contrast to this, to apply Theorem 5.7 we should know from the very beginning that
the sets of vectors  $n_1,\dots,n_m$ and numbers $f_1,\dots,f_m$
correspond to some really existing polyhedral herisson in  ${\bf R}^3$.
Only after this, Theorem 5.7 shows to which extent we can vary the areas of faces
of the polyhedral herisson.

{\bf 5.9.} The {\bf proof} of Theorem 5.7 is obtained below by methods
close to those used in \cite{Al} for proving Theorem 2.8.

By definition, put  $g_j(t)=(1-t)f_j+tg_j$ for $0\leq t\leq 1$ and $1\leq j\leq m$.
Let $T$ be the subset of the segment $[0,1]\subset{\bf R}$
that consists of real numbers $t$ for each of which there exists a
polyhedral herisson in ${\bf R}^3$ which is parallel to and of the same orientation
as $(P; n_1,\dots,n_m)$ and for which the real numbers
$g_1(t),\dots,g_m(t)$ are the oriented areas of the faces.

Show that $T$ is nonempty, closed, and open.
Clearly, in this case $T$ coincides with $[0,1]$ and, in particular,
contains the point  $t=1$.
Obviously, the latter statement is  equivalent to the conclusion of Theorem 5.7.

Trivially, $T$ contains the point $t=0$ and thus is nonempty.

Now prove that $T$ is closed.
Suppose a sequence of points $t_1,\dots,t_k,\dots$ is contained in $T$
and converges to some point $t_0$.
It is sufficient to prove that $t_0\in T$.
Since $t_k\in T$ there exists a polyhedral herisson $(P_k; n_1,\dots,n_m)$ in
${\bf R}^3$ which is parallel to and of the same orientation as
$(P; n_1,\dots,n_m)$ and for which the numbers $g_1(t),\dots,g_m(t)$
are the oriented areas of the faces.

Let $v_1$ be a vertex of $(P_1; n_1,\dots,n_m)$.
Denote by $v_k$, $k\in{\bf N}$, the vertex of $(P_k; n_1,\dots,n_m)$
that has the same spherical image as $v_1$.
Denote by $(\widetilde{P_k}; n_1,\dots,n_m)$ the image of $(P_k; n_1,\dots,n_m)$
under a translation that superposes $v_k$ and $v_1$.

Note that the numbers $g_1(t_k),\dots,g_m(t_k)$ are uniformly bounded from
above for $k\in{\bf N}$.
Besides, $(P; n_1,\dots,n_m)$ is in a general position.
Hence none of the faces of $(P_k; n_1,\dots,n_m)$, $k\in{\bf N}$,
has mutually parallel sides.
It follows from Lemma 5.3 that the perimeters of all faces of all
$(P_k; n_1,\dots,n_m)$, $k\in{\bf N}$, are uniformly bounded from above.
Thus all $(\widetilde{P_k}; n_1,\dots,n_m)$, $k\in{\bf N}$, are contained
in a ball in ${\bf R}^3$.

Fixing each of the vectors $n_1,\dots,n_m$ in turn and applying
Blaschke's selection theorem \cite{Bl} to the sequence of faces
(i.e. convex polygons) of $(\widetilde{P_k}; n_1,\dots,n_m)$, $k\in{\bf N}$,
equipped with the above-chosen normal vector, we conclude that there exists
a subsequence of polyhedral herissons
$(\widetilde{P_{k_j}}; n_1,\dots,n_m)$, $j\in{\bf N}$, such that
every sequence of faces equipped with a common normal vector converges.
It is clear that the limit of each of these sequences of faces is a convex polygon
and the limit of the sequence $(\widetilde{P_{k_j}}; n_1,\dots,n_m)$, $j\in{\bf N}$,
is a polyhedral herisson $(P_0; n_1,\dots,n_m)$, for which the numbers
$g_1(t_0),\dots,g_m(t_0)$ are the oriented areas of the faces.
This yields $t_0\in T$ and thus $T$ is closed.

Finally, prove that $T$ is open.
We need the following auxiliary constructions.

Each polyhedral herisson parallel to a given polyhedral herisson
$(P; n_1,\dots,n_m)$ can be characterized by its support numbers $h_1,\dots,h_m$.
It follows that the set of all polyhedral herissons parallel to
$(P; n_1,\dots,n_m)$ can be considered as a subset of ${\bf R}^m$.
This subset is open because small displacements of planes of the faces
do not destroy the faces.
Unite all the parallel and congruent polyhedral herissons as a
single object called a class.
Since a translation in ${\bf R}^3$ is determined by three parameters,
each class is determined by $m-3$ parameters.
Thus the set of all classes can be treated as an $(m-3)$-dimensional manifold
denoted by ${\bf A}$.

Denote by ${\bf B}$ the set of all $m$-tuples of nonzero real numbers
$f_1,\dots,f_m$ satisfying the condition $\sum_{j=1}^{m} f_jn_j=0$,
where the vectors $n_1,\dots,n_m$ are supposed to be fixed.
The set ${\bf B}$ can be considered as a subset of ${\bf R}^m$, namely as
the $(m-3)$-dimensional plane $\sum_{j=1}^{m} f_jn_j=0$ from which
the points are deleted that have at least one zero coordinate.

Define a mapping $\varphi$ on ${\bf A}$ that takes each polyhedral herisson
parallel to and of the same orientation as $(P; n_1,\dots,n_m)$
(or the set of its support numbers $h_1,\dots,h_m$) to
the set  $f_1,\dots,f_m$  of the oriented areas of its faces.
Lemma 5.1 implies that $\varphi$ maps  ${\bf A}$ into ${\bf B}$.
We have seen above that ${\bf A}$ and ${\bf B}$ are of the same dimensions.
Since the areas of faces depend continuously on location of planes containing
the faces and thus on the support numbers, it follows that $\varphi$ is continuous.
Uniqueness Theorem 4.6 implies that $\varphi$ is injective.

Recall the well-known Brouwer domain invariance theorem:
If $U$ is an open set in ${\bf R}^k$ and $\psi : U\to{\bf R}^k$ is a continuous
one-to-one mapping, then $\psi (U)$ is an open subset of ${\bf R}^k$.
Applying this theorem to $\varphi$ we conclude that $\varphi ({\bf A})$ is an open
subset of ${\bf B}$.

Interpreting the set of numbers $g_1(t),\dots,g_m(t)$ as a point of ${\bf B}$,
we see that, as $t\to t_0$, these points tend to the point
$g_1(t_0),\dots,g_m(t_0)$ of $\varphi ({\bf A})$.
Since $\varphi ({\bf A})$ is open, this implies that, for all $t\in [0,1]$
sufficiently close to $t_0$, the point $g_1(t),\dots,g_m(t)$ lies in $\varphi ({\bf A})$.
Thus all $t\in [0,1]$ sufficiently close to $t_0$ lie in $T$, i.e., $T$ is open.
$\framebox{\phantom{o}}$

{\bf 5.10. Example.}
We show that Theorem 5.7 cannot be directly extended to polyhedral herissons
which are not in a general position.

The nonconvex 11-hedron presented in Fig.~7 is the union of two regular tetrahedra
and a regular triangular prism (all the three are treated as bodies here)
some of whose axes and planes of symmetry mutually coincide.
Denote by $P$ the boundary of this 11-hedron that is a nonconvex
sphere-homeomorphic surface.
The faces of $P$ that do not intersect the triangular prism are called
the {\it bases} of $P$ while the common faces of $P$ and the triangular prism
are called the {\it waist} of $P$.

$$
\epsfbox{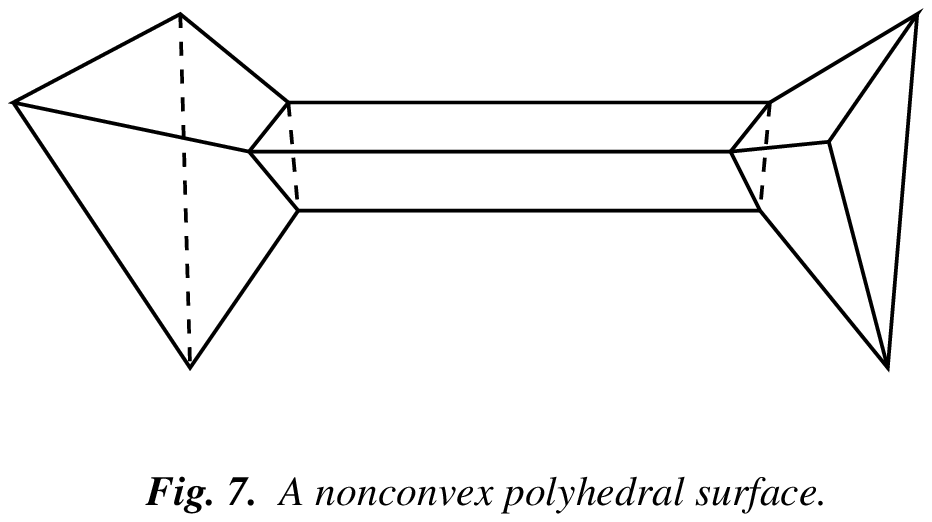}
$$

Make $P$ a herisson $H=(P;n_1,\dots,n_{11})$ by assigning the unit
outward normal vectors $n_1$ and $n_{11}$ to the bases of $P$ and the unit
inward normal vectors  $n_2,\dots,n_{10}$ to all other faces.
Orient $P$ by unit outward normal vectors.
Then the bases of $P$ are its positive faces while all the other faces are
negative.

To be more precise, assume that $P$ is constructed from regular
tetrahedra with unit edge length and the area of the waist is equal to 1.

To each of the eleven vectors $n_j$, $j=1,\dots,11$, assign a number
$g_j$, $j=1,\dots,11$, according to the following rule:
assign the number $\sqrt{3}/4$ (i.e. the area of an isosceles triangle
with side length 1) both to $n_1$ and $n_{11}$ (i.e. to the vectors
attached to the bases of $P$);
assign the number $-1/3$ to each of those vectors $n_j$ that correspond
to the waist of $P$;
and assign the number $-\sqrt{3}/4$ to the rest vectors $n_j$.

We treat the numbers $g_1,\dots,g_{11}$ as the limit values
of the oriented areas of the faces of polyhedral herissons
$H_k=(P_k;n_1,\dots,n_{11})$ constructed as follows.
Given $k\in{\bf N}$, $P_k$ is combinatorially equivalent to $P$,
the equipment of $P_k$ coincides with that of $P$,
the length of an orthogonal cross-section of the waist of $P_k$ equals $1/k$,
and the length of the waist of $P_k$ is equal to $k$.
The oriented areas of the rest faces of $P_k$ either equal $\sqrt{3}/4$
or tend to  $-\sqrt{3}/4$ as $k\to\infty$.
This implies that conditions (1) and (2) of Theorem 5.7 hold true for
the set of numbers $g_1,\dots,g_{11}$.

However prove that there is no polyhedral herisson
that is parallel to and of the same orientation as $H$
and is such that the numbers $g_1,\dots,g_{11}$
are the oriented areas of its faces.

Assume the contrary. Then there exists a polyhedral herisson $H_0$ that is parallel to
and of the same orientation as $H$ and is such that the numbers $g_1,\dots,g_{11}$
are the oriented areas of its faces.
Using the notations introduced in Proof 5.9 of Theorem 5.7,
denote by $H_{\bf A}$ the point of ${\bf A}$ that corresponds to the polyhedral herisson
$H_0$ (i.e. the set of its support numbers) and denote by $H_{\bf B}$ the point of
${\bf B}$ that corresponds to $H_0$ (i.e. the set of the oriented areas of the faces
of $H_0$).
Let $U$ be a neighborhood of $H_{\bf A}$ in ${\bf A}$ such that the sum of
the absolute values of the support numbers of any polyhedral herisson from $U$
is less than twice the sum of the absolute values of the support numbers of $H_0$.
From Proof 5.9 it follows that $\varphi (U)$ is an open neighborhood of the point
$H_{\bf B}$.
In particular, $\varphi (U)$ contains infinitely many polyhedral herissons
(more precisely, the corresponding tuples of the oriented areas of their faces) from
the sequence $H_k$, $k\in{\bf N}$.
Uniqueness Theorem 4.6 implies that only one point of the manifold ${\bf A}$
(or, in terms of Proof 5.9, only one class) corresponds to the polyhedral
herisson $H_k\in{\bf B}$ and, obviously, that a single point lies in $U$.
Consequently there is an infinite subsequence of $H_k$, $k\in{\bf N}$,
with the sum of the absolute values of the support numbers being uniformly bounded
from above.
The latter contradicts the fact that the length of the waist of $H_k$ equals $k$.
This contradiction shows that a polyhedral herisson $H_0$
with the above-indicated properties does not exist.

\section{Concluding remarks}

In this section we briefly discuss several problems on polyhedral herissons
that are beyond the scope of the paper.

{\bf 6.1. Remarks.}
Recall that a convex polytope is called a {\it parallelotope} if the
3-dimensional space ${\bf R}^3$ can be filled with its parallel copies in such a way
that every two copies either do not intersect, or
have a common vertex, or a common edge (a whole edge rather than a part of an edge),
or a common face (an entire face again) \cite{Al}.
Parallelotopes play an important role in such different fields of sciences as
number theory, crystallography, stochastic geometry, theory of vector measures, etc.
We only mention that H.~Minkowski proved his famous uniqueness theorem
(see Theorem 2.6 above) with a view to prove that every parallelotope has a center
of symmetry.

There are a complete classification of parallelotopes in ${\bf R}^3$
(see, for example, \cite{Al}) and several isolated attempts to study nonconvex
space-filling bodies (see, for example, \cite{SS}).
It turns out that there are space-filling nonconvex polyhedral herissons.
A simple example can be constructed as follows.
Consider an equilateral trapezium (as a convex body in the plane)
such that the angle between the bigger base and each of the lateral sides equals $\pi/4$,
and the smaller base is 3 times shorter than the bigger base.
Denote by $P$ the nonconvex body that is the union of the trapezium
and its image under the symmetry with respect to the straight line passing through
the smaller base of the trapezium (Fig.~8).
Denote by $Q$ the Cartesian product of $P$ with a straight line segment
not contained in the plane passing through $P$.
The reader will have no difficulty in showing that the boundary of $Q$
can be made a polyhedral herisson and parallel copies of the herisson
can fill the space.

$$
\epsfbox{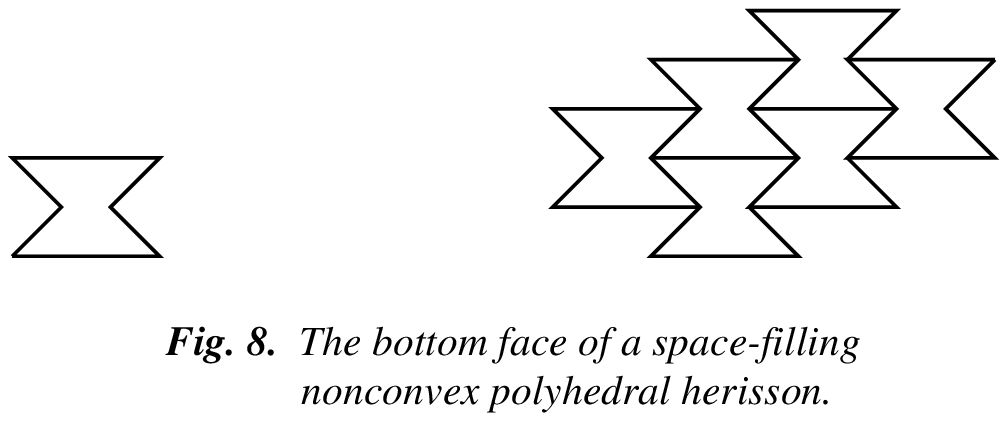}
$$
\medskip

{\bf 6.2. Problem.}
Does there exist a classification of space-filling polyhedral herissons in
${\bf R}^3$?

{\bf 6.3. Remarks.}
To prove his famous uniqueness theorem (see Theorem 2.6 above), H.~Minkowski
developed a machinery based on an operation that is now known as the Minkowski
sum of convex bodies and made an enormous advance in the theory of isoperimetric
inequality by proving, for example, the so-called Minkowski inequality for mixed
volumes (see, for example, \cite{Al}, \cite{Bl}, or \cite{Sc}).

Let $P$ and $Q$ be convex polytopes and let $R$ be their Minkowski sum.
It is well known that $R$ is a convex polytope too.
Moreover, a face of $R$ with a prescribed outward normal vector $n$
can be obtained (up to translation) as the Minkowski sum of the faces
(which may degenerate into edges or vertices) of $P$ and $Q$ with outward
normal vector $n$.

This property leads naturally to a definition of the Minkowski sum of
nonconvex polyhedral herissons.
Actually, faces of polyhedral herissons are convex polygons
and the Minkowski sum is already defined for them.
On the other hand, from the definition of a polyhedral herisson it follows
that the notion of a face with a given outward normal vector is defined
properly.

{\bf 6.4. Problem.}
For polyhedral herissons, do there exist inequalities similar to
those for mixed volumes of convex bodies?

{\bf 6.5. Remarks.}
In \cite{Pa}, the notion of the mixed volume is introduced for virtual polytopes.
As we have mentioned above, the notion of a virtual polytope generalizes that
of a polyhedral herisson.
Nevertheless, no inequality for mixed volumes is given in \cite{Pa}.
Isoperimetric-type inequalities for ``smooth'' herissons are established in \cite{Ma4}.

{\bf 6.6. Problem.}
Does there exist a many-dimensional generalization of the notion of a
polyhedral herisson in ${\bf R}^3$ that allows us to prove a rigidity
theorem similar to theorems proven in \cite{RR}
and to show that for every cross-polytope in ${\bf R}^5$, there
is an equipment which makes it a herisson?

{\bf 6.7. Remark.}
An affirmative answer to Problem 6.6 will imply that there are no flexible
cross-polytopes in the Euclidean 5-space.
This result will support the hypothesis that there are no flexible polyhedra
in ${\bf R}^n$ for $n\geq 5$.
Note that the largest collection of examples of flexible cross-polytopes in
the Euclidean 4-space may be found in \cite{Sta}.

\end{document}